\documentclass[letterpaper, 10 pt, journal, twoside]{IEEEtran}  

\usepackage{amsthm}
\usepackage{amsmath}               
{
  \theoremstyle{plain}
  \newtheorem{assumption}{Assumption}
  \newtheorem{definition}{Definition}
  \newtheorem{theorem}{Theorem}
  
  \newtheorem{lemma}{Lemma}
  
}
\usepackage{amssymb}
\usepackage{bm}
\usepackage{ifthen}
\usepackage{mathtools}
\usepackage{verbatim}
\usepackage{graphicx}
\usepackage{tabularx}
\usepackage[]{hyperref}
\usepackage{xcolor}
\usepackage{algorithm, algorithmicx}
\usepackage[noend]{algpseudocode}
\usepackage{pgffor}
\usepackage{subcaption}
\usepackage{etoolbox}
\usepackage{multirow}
\usepackage{xspace}
\usepackage{refcount}

\algrenewcommand\algorithmicindent{0.9em}%
\algnewcommand{\algorithmicgoto}{\textbf{go to} line}%
\algnewcommand{\Goto}[1]{\algorithmicgoto~\ref{#1}}%

\newcolumntype{C}{>{\centering\arraybackslash}X}
\newcolumntype{L}{>{\raggedright\arraybackslash}X}
\newcolumntype{R}{>{\raggedleft\arraybackslash}X}

\usepackage{xcolor}
\definecolor{commentgreen}{RGB}{0, 153, 0}

\newcommand{\eccalg}{\cite[Algorithm~2]{Malyuta2019a}\xspace}

\newcommand{\BEAS}{\begin{eqnarray*}}
\newcommand{\EEAS}{\end{eqnarray*}}
\newcommand{\BEA}{\begin{eqnarray}}
\newcommand{\EEA}{\end{eqnarray}}
\newcommand{\BEQ}{\begin{equation}}
\newcommand{\EEQ}{\end{equation}}
\newcommand{\BIT}{\begin{itemize}}
\newcommand{\EIT}{\end{itemize}}
\newcommand{\BNUM}{\begin{enumerate}}
\newcommand{\ENUM}{\end{enumerate}}

\newcommand{\BA}{\begin{array}}
\newcommand{\EA}{\end{array}}

\newcommand{\reals}{\mathbb R}
\newcommand{\ball}{\mathbb B}

\newcommand{\binary}{\mathbb I}

\newcommand{\Co}{{\mathop \mathrm{co}}}

\newcommand{\find}{\mathop{\rm find}}

{\begin{quote}}{\end{quote}}

\newcounter{oursection}

\makeatletter
\def\optimize{\@ifnextchar[{\@with}{\@without}}
\def\@with[#1]#2#3#4#5{
  \begin{aligned}
    #2 =~ & #3_{#4}~#5~\mathrm{s.t.} \\
         & #1
  \end{aligned}  
}
\def\@without#1#2#3#4{
  \begin{equation}
    #1 = #2_{#3}~#4
  \end{equation}
}
\makeatother

\definecolor{darkolivegreen}{rgb}{0.33, 0.42, 0.18}

\newcommand{\Behcet}{Beh\c{c}et}
\newcommand{\Acikmese}{A\c{c}{\i}kme\c{s}e}

\newcommand{\transp}{{\scriptscriptstyle\mathsf{T}}}

\newcommand{\kind}{short} 
\newcommand{\colored}{false} 
\newcommand{\longtxt}[1]{%
  \ifthenelse{\equal{\kind}{long} \OR \equal{\kind}{both}}{%
    \ifthenelse{\equal{\colored}{true}}{%
      {\color{red} #1}%
    }{%
      #1%
    }%
  }{}%
}
\newcommand{\shorttxt}[1]{%
  \ifthenelse{\equal{\kind}{short} \OR \equal{\kind}{both}}{%
    \ifthenelse{\equal{\colored}{true}}{%
      {\color{blue} #1}%
    }{%
      #1%
    }%
  }{}%
}

\usepackage{graphics} 

\graphicspath{{./figures/}}

\title{Approximate Multiparametric Mixed-integer Convex Programming}

\author{Danylo Malyuta$^{1}$ and \Behcet{} \Acikmese{}$^{2}$%
  \thanks{$^{1}$Ph.D. student, W.E. Boeing Department of Aeronautics \&
    Astronautics, University of Washington, Seattle, WA 98195, USA
    \texttt{danylo@uw.edu}}%
  \thanks{$^{2}$Professor, W.E. Boeing Department of Aeronautics \&
    Astronautics, University of Washington, Seattle, WA 98195, USA
    \texttt{behcet@uw.edu}}%
}

\begin{document}

\maketitle
\thispagestyle{empty}
\pagestyle{empty}

\begin{abstract}
  We propose an algorithm for generating explicit solutions of multiparametric
  mixed-integer convex programs to within a given suboptimality tolerance. The
  algorithm is applicable to a very general class of optimization problems, but
  is most useful for hybrid model predictive control, where on-line
  implementation is hampered by the worst-case exponential complexity of
  mixed-integer solvers. The output is a simplicial partition which defines a
  static map from the current state to a suboptimal solution. The primary
  theoretical contribution of this paper is to introduce a non-zero optimal cost
  overlap metric which is necessary and sufficient for convergence. The overlap
  size is also linked to partition complexity. The algorithm is massively
  parallelizable and our implementation, which is publicly available, is run on
  a cluster of several hundred processors. Not only does our solution have a
  deterministic runtime, simulations show that our approach is faster than
  on-line optimization by up to three orders of magnitude.
\end{abstract}

\section{Introduction}
\label{introduction}

\IEEEPARstart{H}{ybrid} model predictive control (MPC) handles systems with
discrete switches or piecewise affinely approximated nonlinearities like
chemical powerplants, pipelines and aerospace vehicles
\cite{Bemporad1999a,Blackmore2012,Schouwenaars2006}. This requires solving a
mixed-integer convex program (MICP), which is hampered by the worst-case
exponential complexity of mixed-integer solvers. In this paper, we present a
provably convergent algorithm for computing explicit solutions of MICPs with a
specified suboptimality tolerance.


Several approaches have been proposed to improve hybrid MPC performance. By
leveraging the polynomial runtime complexity of convex solvers, successive
convexification is able to solve nonlinear programs in real-time
\cite{Mao2017}. Recently, the method was extended to handle binary decision
making via state-triggered constraints \cite{Szmuk2019b}. Hence, at least some
MICPs are solvable in real-time. However, this is a local method which may not
always converge to a feasible solution.

The traditional method of ensuring real-time MPC performance while guaranteeing
convergence and global optimality has been to pre-compute the optimal solution
off-line. Various explicit MPC methodologies have been proposed
\cite{Alessio2009}. For MPC laws more complicated than linear or quadratic
programs, exact explicit solutions are generally not possible due to
non-convexity of common active constraint sets \cite{Bemporad2006b}. Instead,
approximate solutions have been proposed via local linearization
\cite{Pistikopoulos2007a} or via optimal cost bounding by affine functions over
simplices \cite{Bemporad2006b} and hyperrectangles \cite{Johansen2004}. An
approximate explicit solution to mixed-integer quadratic programs has been
proposed based on difference-of-convex programming \cite{Alessio2006} and for
MICPs based on local linearization and primal/master subproblems \cite{Dua1998}.

Our contribution in this paper is twofold. First, we introduce a massively
parallelizable algorithm for computing the explicit solution to a very general
class of multiparametric MICPs and to within a user-specified suboptimality
tolerance. Our implementation is available online (see
Section~\ref{sec:simulations}). Second, we define a novel cost overlap metric
and show that it is both the fundamental driver of partition complexity and the
quantity whose non-zero value is necessary and sufficient for convergence. To
the best of our knowledge, the overlap metric is the best theoretical insight
to-date about explicit MPC partition complexity. Furthermore, parallelization of
explicit MPC algorithms has not yet been exploited in existing literature.

The paper is organized as follows. Section~\ref{sec:problem_formulation} defines
the class of programs that our algorithm can handle. Section~\ref{sec:phase2}
presents the solution algorithm. Section~\ref{sec:properties} proves its
convergence and complexity properties. Section~\ref{sec:simulations} applies the
method to the robust control of a satellite's
position. Section~\ref{sec:conclusion} concludes with future research
directions.

\textit{Notation}: $\binary\triangleq\{0,1\}$ is the binary set and
$\ball\triangleq\{x:\|x\|_2\le 1\}$ is the unit ball. Matrices are uppercase
(e.g. $A$), scalars, vectors and functions are lowercase (e.g. $x$), and sets
are calligraphic uppercase (e.g. $\mathcal S$). $\Co\mathcal S$,
$\mathcal S^{\mathrm{c}}$, $\partial\mathcal S$ and $\mathcal V(\mathcal R)$
denote respectively the convex hull, complement, boundary and extreme points
(e.g. vertices) of $\mathcal S$. The cardinality of a countable set $\mathcal S$
is $|\mathcal S|$. Given $\mathcal A\subseteq\reals^n$, $b\in\reals^n$ and
$s\in\reals$, $\mathcal A+b\triangleq\{a+b\in\reals^n:a\in\mathcal A\}$ and
$s\mathcal A\triangleq\{sa: a\in\mathcal A\}$.

\section{Problem Formulation}
\label{sec:problem_formulation}

Our algorithm can handle any problem that can be formulated as the following
multiparametric MICP:
\begin{equation}
  \label{eq:minlp}
  \optimize[
  g(\theta,x,\delta)=0,~h(\theta,x,\delta)\in\mathcal K,
  ]{V^*(\theta)}{\min}{x,\delta}{f(\theta,x,\delta)}
  \tag{P$_\theta$}
\end{equation}
where $\theta\in\reals^p$ is a parameter, $x\in\reals^n$ is a decision vector
and $\delta\in\binary^m$ is a binary \textit{commutation}. The cost function
$f:\reals^p\times\reals^n\times\reals^m\to\reals$ is jointly convex and the
constraint functions $g:\reals^p\times\reals^n\times\reals^m\to\reals^{l}$ and
$h:\reals^p\times\reals^n\times\reals^m\to\reals^d$ are affine in their first
two arguments. The convex cone
$\mathcal K=\mathcal C_1\times\cdots\times\mathcal C_{q}\subset\reals^d$ is a
Cartesian product of $q$ convex cones. Examples include the positive orthant,
the second-order cone and the positive semidefinite cone. If $\delta$ is fixed,
\eqref{eq:minlp} becomes a multiparametric convex program:
\begin{equation}
  \label{eq:nlp}
  \optimize[
  g(\theta,x,\delta)=0,~h(\theta,x,\delta)\in\mathcal K.
  ]{V^*_\delta(\theta)}{\min}{x}{f(\theta,x,\delta)}
  \tag{P$_{\theta}^\delta$}
\end{equation}

Let $\Theta^*\subseteq\reals^p$ and $\Theta_\delta^*\subseteq\Theta^*$ denote
respectively the parameter sets for which \eqref{eq:minlp} and \eqref{eq:nlp}
are feasible. Define the following three maps similarly to
\cite{Bemporad2006b,Malyuta2019a}.

\begin{definition}
  The optimal map $f_\delta^*:\Theta^*\to\binary^m$ associates
  $\theta\in\Theta^*$ to any optimal commutation of (\ref{eq:minlp}), that is
  any $\delta\in\{\delta\in\binary^m:V^*(\theta)=V^*_\delta(\theta)\}$.
\end{definition}

\begin{definition}
  The feasible map $f_\delta:\Theta^*\to\binary^m$ associates
  $\theta\in\Theta^*$ to a commutation such that (\ref{eq:nlp}) is feasible.
\end{definition}

\begin{definition}
  \label{definition:suboptimal_commutation_function}
  The suboptimal map $f_\delta^{\epsilon}:\Theta^*\to\binary^m$ associates
  $\theta\in\Theta^*$ to an $\epsilon$-suboptimal commutation $\delta$ such that
  \begin{equation}
    \label{eq:epsilon_suboptimality}
    V^*_\delta(\theta)-V^*(\theta)<\max\{\epsilon_{\mathrm{a}},
    \epsilon_{\mathrm{r}}V^*(\theta)\},
  \end{equation}
  where $\epsilon_{\mathrm{a}}$ and $\epsilon_{\mathrm{r}}$ are the absolute and
  relative errors.
\end{definition}

The next section presents an algorithm for computing $f_{\delta}^\epsilon$ over
a subset $\Theta\subseteq\Theta^*$. It is assumed that $\Theta$ is a
full-dimensional convex polytope in vertex representation. One can choose
$\Theta$ following the advice of \cite[Section~IV-C]{Malyuta2019a}.

\section{Explicit Solution of \eqref{eq:minlp}}
\label{sec:phase2}

\subsection{Suboptimal Map Computation}

We begin by computing $f_\delta^\epsilon$ as a simplicial partition
$\mathcal P\triangleq\{(\mathcal R_i,\delta_i,\{x_{i,j}^*\}_{j=1}^{|\mathcal V(\mathcal R_i)|})\}_{i=1}^{P}$
of $\Theta$. Each simplex $\mathcal R_i$ is associated with an
$\epsilon$-suboptimal $\delta_i$ and the set
$\{x_{i,j}^*\}_{j=1}^{|\mathcal V(\mathcal R_i)|}$ of optimal decision vectors
of \eqref{eq:nlp} at the vertices of $\mathcal R_i$.

Algorithm~\ref{alg:lcss} stores $f_\delta^\epsilon$ as a binary tree and
computes it as follows. To initialize, all leaves of the binary tree $f_\delta$
output by \eccalg are converted into nodes, i.e. elements that have children in
the final tree. Lines
\ref{alg:lcss:line:check_suboptimality}-\ref{alg:lcss:line:parwhile_end} carry
out the main work of partitioning the simplex-commutation tuple
$(\mathcal R,\delta)$. First, the algorithm checks if $\delta$ is
$\epsilon$-suboptimal in $\mathcal R$. If it is not, the following mixed-integer
program must be feasible due to \eqref{eq:epsilon_suboptimality}:
\begin{equation}
  \label{eq:suboptimality_check_nonconvex}
  \optimize[
  V_{\delta}^*(\theta)-V_{\delta'}^*(\theta)\ge
  \max\{\epsilon_{\mathrm{a}},
  \epsilon_{\mathrm{r}}V_{\delta'}^*(\theta)\}.
  ]{\delta^*,\theta^*}{\find}{\theta\in\mathcal R}{\delta'}
  \tag{$\textnormal{E}_{\delta}^{\mathcal R}$}
\end{equation}
The costs $V_\delta^*$ and $V_{\delta'}^*$ are convex
\cite[Lemma~1]{Malyuta2019a}. However, $V_\delta^*$ appears on the wrong side of
the inequality. Hence, \eqref{eq:suboptimality_check_nonconvex} is non-convex
and so is not readily solvable. As a remedy, we formulate a conservative convex
upper bound.

\begin{definition}
  \label{definition:cost_affine_over_approximator}
  Let $v_i\in\mathcal V(\mathcal R)$ be the $i$-th vertex of $\mathcal R$ and
  let $\theta=\sum_{i=1}^{|\mathcal V(\mathcal R)|}\alpha_iv_i$ where
  $\alpha_i\ge 0$ and $\sum_{i=1}^{|\mathcal V(\mathcal R)|}\alpha_i=1$. 
  The affine over-approximator of $V^*_\delta(\theta)$ over $\mathcal R$ is:
  \begin{equation}
    \label{eq:over_approximator}
    \bar V_\delta(\theta)\triangleq\sum_{i=1}^{|\mathcal V(\mathcal R)|}\alpha_iV^*_\delta(v_i).
  \end{equation}
\end{definition}

Since $V_\delta^*$ is convex, $V_\delta^*(\theta)\le\bar V_\delta(\theta)$
$\forall\theta\in\mathcal R$. Hence, the following problem is convex and
``conservative'' in the sense of
Theorem~\ref{theorem:epsilon_suboptimality_sufficiency}:
\begin{equation}
  \label{eq:suboptimality_check}
  \optimize[
  \bar V_{\delta}(\theta)-V_{\delta'}^*(\theta)\ge
  \max\{\epsilon_{\mathrm{a}},
  \epsilon_{\mathrm{r}}V_{\delta'}^*(\theta)\}.
  ]{\delta^*,\theta^*}{\find}{\theta\in\mathcal R}{\delta'}
  \tag{$\bar{\textnormal{E}}_{\delta}^{\mathcal R}$}
\end{equation}

\begin{theorem}
  \label{theorem:epsilon_suboptimality_sufficiency}
  If \eqref{eq:suboptimality_check} is infeasible then $\delta$ is
  $\epsilon$-suboptimal.

  \begin{proof}
    If $\delta$ is not $\epsilon$-suboptimal then
    \eqref{eq:suboptimality_check_nonconvex} is feasible, hence
    \eqref{eq:suboptimality_check} is feasible. By contraposition, if
    \eqref{eq:suboptimality_check} is infeasible then $\delta$ must be
    $\epsilon$-suboptimal.
  \end{proof}
\end{theorem}

Since \eqref{eq:suboptimality_check} is a MICP, its feasibility can be certified
with a mixed-integer solver. If \eqref{eq:suboptimality_check} is infeasible,
line~\ref{alg:lcss:line:close} converts the node to a leaf since no further
partitioning is necessary. Otherwise, we cannot conclude about the
$\epsilon$-suboptimality of $\delta$. In this case, we search for a potentially
better commutation $\delta^*$ via the following extension of
\eqref{eq:suboptimality_check}:
\begin{equation}
  \label{eq:better_delta_selection}
  \optimize[
  \bar V_{\delta}(\theta)-V_{\delta'}^*(\theta)\ge
  \max\{\epsilon_{\mathrm{a}},
  \epsilon_{\mathrm{r}}V_{\delta'}^*(\theta)\}, \\
  &\delta'\in\{\delta''\in\binary^m\setminus\{\delta\}:\mathcal R\subseteq\Theta_{\delta''}^*\}.
  ]{\delta^*,\theta^*}{\find}{\theta\in\mathcal R}{\delta'}
  \tag{$\bar{\textnormal{D}}_{\delta}^{\mathcal R}$},
\end{equation}
where the last constraint ensures that $\delta'$ is feasible in $\mathcal R$ and
can be embedded via \cite[Lemma~2]{Malyuta2019a}. If
\eqref{eq:better_delta_selection} is infeasible, Section~\ref{sec:properties}
shows that a sound strategy is to keep $\delta$ and to split $\mathcal R$ in
half at the midpoint of its longest edge, as done on
line~\ref{alg:lcss:line:split_end}. This may also be done if
\eqref{eq:better_delta_selection} is feasible since shrinking $\mathcal R$
improves the accuracy of the over-approximator
\eqref{eq:over_approximator}. However, to avoid unnecessary partitioning,
line~\ref{alg:lcss:line:check_variability} checks if $V_\delta^*$ varies over
$\mathcal R$ by less than the $\epsilon$-suboptimality threshold with respect to
$V_{\delta^*}(\theta^*)$ as output by \eqref{eq:better_delta_selection}:
\begin{equation}
  \label{eq:curvature_constraint}
  \max_{\theta\in\mathcal R}V_\delta^*(\theta)-
  \min_{\theta\in\mathcal R}V_\delta^*(\theta)<\max\{\epsilon_{\mathrm{a}},
  \epsilon_{\mathrm{r}}V_{\delta^*}^*(\theta^*)\},
\end{equation}
which can be done with convex optimization. If \eqref{eq:curvature_constraint}
holds, Theorem~\ref{theorem:convergence} assures that $\mathcal R$ need not be
subdivided further and $\delta^*$ is assigned directly on
line~\ref{alg:lcss:line:exit_without_split}. If \eqref{eq:curvature_constraint}
does not hold, Section~\ref{sec:properties} shows that a sound strategy is to
split $\mathcal R$ in half at the midpoint of its longest edge on
line~\ref{alg:lcss:line:split_end}.

\begin{algorithm}
  \centering
  \begin{algorithmic}[1]
    \State Run \eccalg and relabel all leaves as nodes
    \While{any nodes exist}
    \label{alg:lcss:line:while}
    \State $(\mathcal R,\delta)\gets\text{the most recently added node}$
    \label{alg:lcss:line:parwhile_start}
    \If{\eqref{eq:suboptimality_check} infeasible}
    \label{alg:lcss:line:check_suboptimality}
    \State Change node to leaf
    $(\mathcal R,\delta,\{x_j^*\}_{j=1}^{|\mathcal V(\mathcal R)|})$
    \label{alg:lcss:line:close}
    \Else
    \State $\delta^*,\theta^*\gets\text{solve
      \eqref{eq:better_delta_selection}}$
    \label{alg:lcss:line:better_delta_selection}
    \If{\eqref{eq:better_delta_selection} feasible
      and \eqref{eq:curvature_constraint} holds}
    \label{alg:lcss:line:check_variability}
    \State Change node to $(\mathcal R,\mathcal \delta^*)$
    \label{alg:lcss:line:exit_without_split}
    \Else
    \State $\delta^*\gets\delta\text{ if \eqref{eq:better_delta_selection} infeasible}$
    \State
    $v_1,v_2\gets\arg\max_{v,v'\in\mathcal V(\mathcal R)}\|v-v'\|_2$
    \label{alg:lcss:line:split_start}
    \State $\mathcal S_i\gets\Co\{(\mathcal V(\mathcal R)\setminus\{v_i\})
    \cup\{(v_1+v_2)/2\}\},~i=1,2$
    \State Add child nodes $(\mathcal S_1,\delta^*)$ and
    $(\mathcal S_2,\delta^*)$
    \label{alg:lcss:line:split_end}
    \label{alg:lcss:line:parwhile_end}
    \EndIf
    \EndIf
    \EndWhile
    \caption{Computation of $f_\delta^{\epsilon}$.}
    \label{alg:lcss}
  \end{algorithmic}
\end{algorithm}

\subsection{Explicit Implementation}

At this point the partition $\mathcal P$ is available. Consider a cell
$(\mathcal R,\delta,\{x_{j}^*\}_{j=1}^{|\mathcal V(\mathcal R)|})$ of
the partition. We first recognize the following property.

\begin{theorem}
  \label{theorem:cvx_comb_decision_vector}
  Suppose that $\theta\in\mathcal R$ and let $v_j\in\mathcal V(\mathcal R)$
  be the $j$-th vertex of $\mathcal R$. One can write
  $\theta=\sum_{j=1}^{|\mathcal V(\mathcal R)|}\alpha_j v_j$ where
  $\alpha_j\ge 0$ and $\sum_{j=1}^{|\mathcal V(\mathcal R_j)|}\alpha_j=1$. Let
  $\hat x\triangleq\sum_{j=1}^{|\mathcal V(\mathcal R)|}\alpha_j x_{j}^*$
  and $\hat V_{\delta}(\theta)\triangleq f(\theta,\hat x,\delta)$. Then $\hat x$
  is feasible for \eqref{eq:nlp} and
  \begin{equation}
    \label{eq:f_le_V}
    \hat V_{\delta}(\theta) \le \bar V_{\delta}(\theta).
  \end{equation}

  \begin{proof}
    For feasibility, exploit that $g(\cdot,\cdot,\delta)$ and
    $h(\cdot,\cdot,\delta)$ in \eqref{eq:nlp} are affine. For
    \eqref{eq:f_le_V}, exploit that $f(\cdot,\cdot,\delta)$ is convex:
    \begin{equation*}
      \hat V_{\delta}(\theta) \le
      \sum_{j=1}^{\mathclap{|\mathcal V(\mathcal R)|}}\alpha_j f(v_j,x_{j}^*,\delta) =
      \sum_{j=1}^{\mathclap{|\mathcal V(\mathcal R)|}}\alpha_j V_{\delta}^*(v_j) =
      \bar V_{\delta}(\theta).\qedhere
    \end{equation*}
  \end{proof}
\end{theorem}

As a result of \eqref{eq:f_le_V}, \eqref{eq:suboptimality_check} continues to be
infeasible if $\bar V_\delta(\theta)$ is substituted by
$\hat V_\delta(\theta)$. Therefore $\hat x$ is $\epsilon$-suboptimal, i.e.
\begin{equation}
  \hat V_\delta(\theta)-V^*(\theta)<\max\{\epsilon_{\mathrm{a}},
  \epsilon_{\mathrm{r}}V^*(\theta)\}.
\end{equation}

Hence, given $\theta\in\mathcal R\subseteq\Theta$, $\hat x$ in
Theorem~\ref{theorem:cvx_comb_decision_vector} gives an explicit solution and is
obtained by querying the partition via:
\begin{equation}
  \label{eq:explicit_implementation}
  \hat x = \sum_{j=1}^{\mathclap{|\mathcal V(\mathcal R)|}}\alpha_j x_{j}^*\text{
    where
  }\theta=\sum_{j=1}^{\mathclap{|\mathcal V(\mathcal R)|}}\alpha_jv_j,~v_j\in\mathcal V(\mathcal R).
\end{equation}

\section{Properties}
\label{sec:properties}

\subsection{Convergence}
\label{subsec:convergence}

This section proves that Algorithm~\ref{alg:lcss} converges if
Assumption~\ref{assumption:positive_overlap} holds. Without loss of generality,
we restrict the discussion to $\Delta$, the set of feasible commutations in
$\Theta$.

\begin{definition}
  \label{definition:overlap}
  The \textnormal{overlap} is the largest $\gamma\ge 0$ such that for each
  $\theta\in\Theta$, $\exists\delta\in\Delta$ which is $\epsilon$-suboptimal in
  $(\gamma\ball+\theta)\setminus\Theta^{\mathrm{c}}$.
\end{definition}

\begin{assumption}
  \label{assumption:positive_overlap}
  The overlap is positive, i.e. $\gamma>0$.
\end{assumption}

\begin{figure}
  \centering
  \includegraphics[width=1\columnwidth]{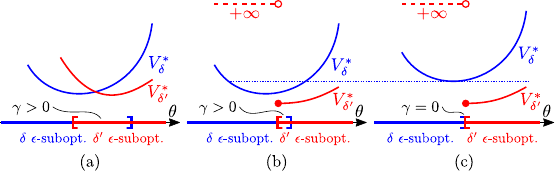}
  \caption{Illustration of ``overlap'' in
    Definition~\ref{definition:overlap}. In (a) the overlap is positive thanks
    to local continuity and in (b) it is positive because the downward jump from
    $V_\delta^*$ to $V_{\delta'}^*$ is not too high. In (c) the jump is too
    high, causing zero overlap.}
  \label{fig:overlap}
\end{figure}

The overlap $\gamma$ is between sets where a given commutation is
$\epsilon$-suboptimal. Its value is a non-trivial property of \eqref{eq:minlp}
which increases for larger values of $\epsilon_{\mathrm{a}}$ and
$\epsilon_{\mathrm{r}}$. Because Algorithm~\ref{alg:lcss} requires \eccalg to
converge, we know that $\gamma\ge 0$ exists. Figure~\ref{fig:overlap}
illustrates just three overlap possibilities. For convergence,
Algorithm~\ref{alg:lcss} should not ``oscillate'' between $\delta$ choices for
the same $\mathcal R$. This is guaranteed by the following lemma.

\begin{lemma}
  \label{lemma:delta_rejection_forever}
  Let $(\mathcal R,\delta)$ be the node selected on
  line~\ref{alg:lcss:line:parwhile_start} at some iteration of
  Algorithm~\ref{alg:lcss}. If $\delta$ is replaced with $\delta^*$ on
  line~\ref{alg:lcss:line:exit_without_split}, the node $(\mathcal R,\delta)$
  will not reappear in a future iteration.

  \begin{proof}
    We begin by showing that
    \begin{equation}
      \label{eq:lemma4_0}
      \min_{\theta\in\mathcal R}V_{\delta^*}^*(\theta)<
      \min_{\theta\in\mathcal R}V_{\delta}^*(\theta).
    \end{equation}

    Since \eqref{eq:better_delta_selection} is feasible,
    $\exists\theta^*\in\mathcal R$ such that
    \begin{equation}
      \label{eq:lemma4_1}
      V_{\delta^*}^*(\theta^*)\le
      \bar V_\delta(\theta^*)-\max\{\epsilon_{\mathrm{a}},
      \epsilon_{\mathrm{r}}V_{\delta^*}^*(\theta^*)\}.
    \end{equation}
    
    Since \eqref{eq:curvature_constraint} holds, we have:
    \begin{align}
      \nonumber
      \bar V_\delta(\theta^*)
      &\le \max_{\theta\in\mathcal R}V_\delta^*(\theta) \\
      \label{eq:lemma4_2}
      &< \min_{\theta\in\mathcal R}V_\delta^*(\theta)+
        \max\{\epsilon_{\mathrm{a}},
        \epsilon_{\mathrm{r}}V_{\delta^*}^*(\theta^*)\}.
    \end{align}

    Substituting \eqref{eq:lemma4_2} into \eqref{eq:lemma4_1} shows that
    \eqref{eq:lemma4_0} holds. If $(\mathcal R,\delta)$ reappears in a future
    iteration, then it must be that all of the preceding iterations finished on
    line~\ref{alg:lcss:line:exit_without_split}. Consider a future iteration
    where the node is $(\mathcal R,\tilde\delta)$. Recursively applying
    \eqref{eq:lemma4_0}, we have:
    \begin{equation}
      \label{eq:lemma4_5}
      \min_{\theta\in\mathcal R} V_{\tilde\delta}^*(\theta)<\cdots
      <\min_{\theta\in\mathcal R} V_{\delta^*}^*(\theta)<
      \min_{\theta\in\mathcal R}V_\delta^*(\theta).
    \end{equation}

    By contradiction, suppose that $\delta$ is chosen on
    line~\ref{alg:lcss:line:better_delta_selection}. This means that
    $\exists\tilde\theta\in\mathcal R$ such that
    \begin{equation}
      \label{eq:lemma4_3}
      V_{\delta}^*(\tilde\theta)\le
      \bar V_{\tilde\delta}(\tilde\theta)-\max\{\epsilon_{\mathrm{a}},
      \epsilon_{\mathrm{r}}V_{\delta}^*(\tilde\theta)\}.
    \end{equation}
    
    In the same way that we obtained \eqref{eq:lemma4_2}, we have:
    \begin{equation}
      \label{eq:lemma4_4}
      \bar V_{\tilde\delta}(\tilde\theta)
      < \min_{\theta\in\mathcal R}V_{\tilde\delta}^*(\theta)+
      \max\{\epsilon_{\mathrm{a}},
      \epsilon_{\mathrm{r}}V_{\delta}^*(\tilde\theta)\}.
    \end{equation}
        
    Using \eqref{eq:lemma4_4} in \eqref{eq:lemma4_3}, we have
    \begin{equation*}
      V_{\delta}^*(\tilde\theta)<
      \min_{\theta\in\mathcal R}V_{\tilde\delta}^*(\theta)
      ~\Rightarrow~
      \min_{\theta\in\mathcal R}V_{\delta}^*(\theta)<
      \min_{\theta\in\mathcal R}V_{\tilde\delta}^*(\theta),
    \end{equation*}
    which contradicts \eqref{eq:lemma4_5}, hence $\delta$ cannot be more optimal
    than $\tilde\delta$ and thus cannot be re-associated with $\mathcal R$.
  \end{proof}
\end{lemma}

\begin{theorem}
  \label{theorem:convergence}
  Algorithm~\ref{alg:lcss} terminates if and only if
  Assumption~\ref{assumption:positive_overlap} holds.

  \begin{proof}
    Suppose that Assumption~\ref{assumption:positive_overlap} holds. The
    algorithm terminates when all nodes become leaves. An iteration can exit on
    line~\ref{alg:lcss:line:close}, \ref{alg:lcss:line:exit_without_split} or
    \ref{alg:lcss:line:split_end}. Since exiting on
    line~\ref{alg:lcss:line:close} terminates a branch, it is necessary and
    sufficient to show that only a finite number of iterations can exit on
    lines~\ref{alg:lcss:line:exit_without_split} or
    \ref{alg:lcss:line:split_end}. Since $V_\delta^*$ is convex, it is
    continuous and therefore \eqref{eq:curvature_constraint} holds for a small
    enough $\mathcal R$ but with a non-empty interior
    \cite{Boyd2004,Royden1988}. If an iteration exits on
    line~\ref{alg:lcss:line:split_end}, the volume of $\mathcal R$ is halved and
    its size is reduced, so after a finite number of iterations $\mathcal R$
    will be small enough such that \eqref{eq:curvature_constraint} holds and
    $\mathcal R\subset\gamma\ball+\theta$ for some $\theta\in\Theta$. Once this
    occurs, by Definition~\ref{definition:overlap} $\exists\delta^*\in\Delta$
    such that \eqref{eq:better_delta_selection} is feasible. As a result, for
    any $\delta\in\Delta$ it will take a finite number of iterations until all
    iterations persistently exit on
    line~\ref{alg:lcss:line:exit_without_split}. However, by
    Lemma~\ref{lemma:delta_rejection_forever} and since $|\Delta|$ is finite,
    this can only occur a finite number of times. Thus, after a finite number of
    iterations there will remain only one possible choice of $\delta$ and the
    iteration will exit on line~\ref{alg:lcss:line:close}, so the algorithm
    terminates. If Assumption~\ref{assumption:positive_overlap} does not hold,
    it will take infinite iterations until
    $\mathcal R\subset\gamma\ball+\theta$, so the algorithm does not terminate.
  \end{proof}
\end{theorem}

Theorem~\ref{theorem:convergence} indicates that the partition complexity is
driven by $\gamma$ and the required ``smallness'' of $\mathcal R$ such that
\eqref{eq:curvature_constraint} holds. We formally define the latter quantity
below.

\begin{definition}
  \label{definition:variability}
  The (conservative) \textnormal{variability} is the largest $\nu>0$ such that
  \eqref{eq:curvature_constraint} holds for any $\delta\in\Delta$, any
  $\mathcal R\subset\nu\ball+\theta$, any $\theta\in\Theta$ and any
  $\theta^*\in\mathcal R$.
\end{definition}

Combining Definitions~\ref{definition:overlap} and \ref{definition:variability},
we can state an overall \textit{condition number}:
\begin{equation}
  \label{eq:condition_number}
  \psi \triangleq \min\{\gamma,\nu\}^{-1},
\end{equation}
which is positively correlated to how much $\mathcal R$ must be subdivided until
Theorem~\ref{theorem:convergence} assures convergence. We call \eqref{eq:minlp}
with small $\psi$ ``well-conditioned'' and Algorithm~\ref{alg:lcss} will
converge faster. Note that larger $\epsilon_{\mathrm{a}}$ and
$\epsilon_{\mathrm{r}}$ decrease $\psi$.

\subsection{Complexity}
\label{subsec:complexity}

This section proves that evaluating \eqref{eq:explicit_implementation} has
polynomial complexity. We assume that $\Theta$ is a simplex, so the partition
$\mathcal P$ is a binary tree. We begin by determining the tree depth.

\begin{lemma}
  \label{lemma:tree_depth_complexity}
  The depth $\tau$ of the tree output by Algorithm~\ref{alg:lcss} is
  $\mathcal O(p^2\log(\psi))$.

  \begin{proof}
    In the worst case, Algorithm~\ref{alg:lcss} has to reduce the size of
    $\mathcal R$ until $\mathcal R\subset\psi\ball+\theta$ for some
    $\theta\in\Theta$. Once this occurs, Theorem~\ref{theorem:convergence}
    assures that a future iteration will close the corresponding branch without
    further subdivision. It was shown in \cite[Theorem~2]{Malyuta2019a} that
    reducing $\mathcal R$ until $\mathcal R\subset\psi\ball+\theta$ takes
    $\tau=\mathcal O(p^2\log(\psi))$ subdivisions.
  \end{proof}
\end{lemma}

\begin{theorem}
  \label{theorem:onlineevalcomplexity}
  The evaluation complexity of \eqref{eq:explicit_implementation} is
  $\mathcal O(p^4)$.

  \begin{proof}
    As explained in Section~\ref{subsec:improving_storage_memory_requirement},
    checking if $\theta\in\mathcal R$ can be done via a matrix-vector product,
    which is $\mathcal O(p^2)$. Since there are $\tau$ such checks to perform
    and since $\tau=\mathcal O(p^2)$ due to
    Lemma~\ref{lemma:tree_depth_complexity}, it takes $\mathcal O(p^4)$
    operations to find the $\mathcal R$ which contains $\theta$. It subsequently
    takes $\mathcal O(p)$ operations to compute $\hat x$, hence the overall
    evaluation complexity is $\mathcal O(p^4)$.
  \end{proof}
\end{theorem}

Theorem~\ref{theorem:onlineevalcomplexity} stands in contrast to implementing
\eqref{eq:minlp} directly with a mixed-integer solver, which has an exponential
runtime $\mathcal O(2^m)$.

\section{Simulation Examples}
\label{sec:simulations}

\begin{figure}
  \centering
  \includegraphics[width=0.5\linewidth]{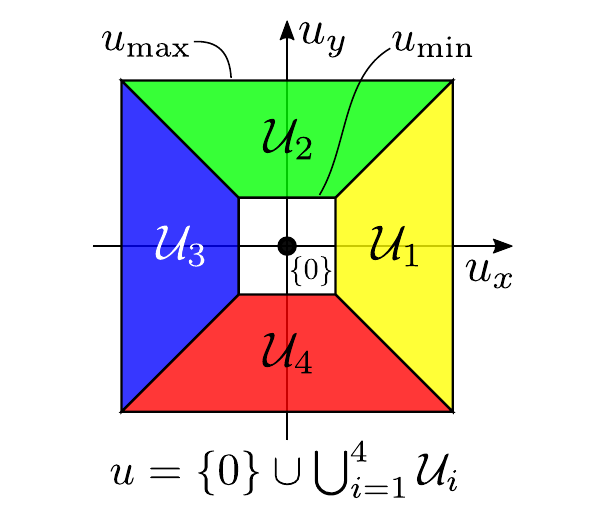}
  \caption{A non-convex input lower-bound constraint for satellite MPC can be
    modeled as the union of convex sets. This induces a mixed-integer program.}
  \label{fig:cwh_union}
\end{figure}

This section presents two examples to corroborate the effectiveness of
Algorithm~\ref{alg:lcss} and the conclusions of Section~\ref{sec:properties}. We
consider robust hybrid MPC of a satellite's out-of-plane (\verb|cwh_z|) and
in-plane (\verb|cwh_xy|) position. Explicit MPC is relevant for satellite
control because the conservative design of space systems typically prohibits the
use of on-line optimization. We use Clohessy-Wiltshire-Hill dynamics:
\begin{align}
  \label{eq:cwh_xy}
  \text{\texttt{cwh\_xy}:~}
  &\left\{
    \begin{array}{@{}ll@{}}
      \ddot x = 3\omega_0^2x+2\omega_o\dot y+u_x+w_x, \\
      \ddot y = -2\omega_o \dot x+u_y+w_u,
    \end{array}\right.
 \\
  \text{\texttt{cwh\_z}:~}
  &\hspace{3.25mm}\ddot z = -\omega_o^2 z+u_z+w_z,
\end{align}
where $\omega_o$ is the orbital rate in~rad/s and $w$ are disturbance terms. The
full explanation is provided in~\cite{Malyuta2019}. For \verb|cwh_xy|, the state
is $(x,\dot x,y,\dot y)\in\reals^4$ and the input is $(u_x,u_y)\in\reals^2$. For
\verb|cwh_z|, the state is $(z,\dot z)\in\reals^2$ and the input is
$u_z\in\reals$. Assuming an impulsive input, the system is discretized at a
$T_{\mathrm{s}}=100$~s thruster firing period. On top of \cite{Malyuta2019}, we
add a lower-bound constraint $\|u\|_\infty\ge u_{\min}$ on the thrust magnitude,
which arises from the thruster impulse-bit. The constraint is non-convex but can
be modeled as the union of convex sets, as illustrated in
Figure~\ref{fig:cwh_union}, yielding a mixed-integer second order cone
program. We use a prediction horizon $N=4$ and, letting
$\mathcal X\triangleq[-10,10]~\text{cm}\times[-1,1]~\text{mm/s}$, choose
$\Theta=\mathcal X$ for \verb|cwh_z| and $\Theta=\mathcal X\times\mathcal X$ for
\verb|cwh_xy|. It was shown in \cite{Malyuta2019} that this $\Theta$ choice is a
robust controlled invariant set, hence the partition will be sufficient for
controlling the satellite.

Our implementation is available
online\footnote{\href{https://github.com/dmalyuta/explicit_hybrid_mpc}{\texttt{https://github.com/dmalyuta/explicit\_hybrid\_mpc}}}
and uses Python 3.7.2, CVXPY 1.0.21 \cite{cvxpy} and MOSEK 9.0.87
\cite{mosek}. Recognizing that
lines~\ref{alg:lcss:line:check_suboptimality}-\ref{alg:lcss:line:parwhile_end}
can run in parallel across tree branches, we used MPICH 3.2 in CentOS 7 on a
cluster of up to 420 2.4~GHz Intel E5-2680 CPU cores with 20~GB of RAM per
compute node (28 cores). The code can also run locally.

\begin{table*}
  \centering
  \begin{tabularx}{0.8\linewidth}{>{\hsize=1.7cm}CCCCCCCC}
    \vspace{1em}
    Example
    & \vspace{1em} $s_{\mathrm{a}}$
    & \vspace{1em} $\epsilon_{\mathrm{r}}$
    & \vspace{1em} $\tau$
    & \vspace{1em} $\lambda$
    & \vspace{1em} $T_{\mathrm{wall}}$ [hr]
    & \vspace{1em} $T_{\mathrm{cpu}}$ [hr]
    & \vspace{1em} $M$ [MB] \\
    \hline \hline
    \verb|cwh_z|  & $0.50$ & $2.00$  & $13$ & $101$     & $0.01$  & $0.09$    & $<0.01$ \\
    \verb|cwh_z|  & $0.25$ & $1.00$  & $17$ & $978$     & $0.06$  & $0.96$    & $<0.01$ \\
    \verb|cwh_z|  & $0.10$ & $0.10$  & $20$ & $13500$   & $0.31$  & $7.72$    & $11$ \\
    \verb|cwh_z|  & $0.03$ & $0.05$ & $26$ & $235231$  & $1.91$  & $154.19$  & $202$ \\
    \verb|cwh_z|  & $0.01$ & $0.01$ & $31$ & $3322941$ & $6.37$  & $2516.98$ & $2916$ \\
    \verb|cwh_xy| & $0.50$ & $2.00$  & $32$ & $30448$   & $0.57$  & $53.44$   & $36$ \\
    \verb|cwh_xy| & $0.25$ & $1.00$  & $49$ & $884323$  & $3.38$  & $1297.35$ & $1069$ \\
    \hline
  \end{tabularx}
  \caption{Numerical results for several $\epsilon$-suboptimality settings;
    $\tau$ is the tree depth, $\lambda$ is the leaf count, $T_{\mathrm{wall}}$
    is the Algorithm~\ref{alg:lcss} runtime, $T_{\mathrm{cpu}}$ is the
    computation time summed across parallel processors, and $M$ is the tree file
    size.}
  \label{tab:results}
\end{table*}

\begin{figure}
  \centering
  \includegraphics[width=0.8\columnwidth]{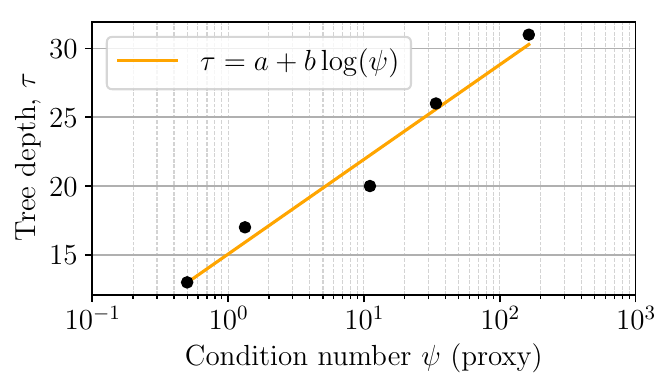}
  \caption{The tree depth for \texttt{cwh\_z} is approximately logarithmic in
    the condition number $\psi$, as predicted by
    Lemma~\ref{lemma:tree_depth_complexity}.}
  \label{fig:depth_complexity}
\end{figure}

\begin{figure}
  \centering
  \includegraphics[width=1\linewidth]{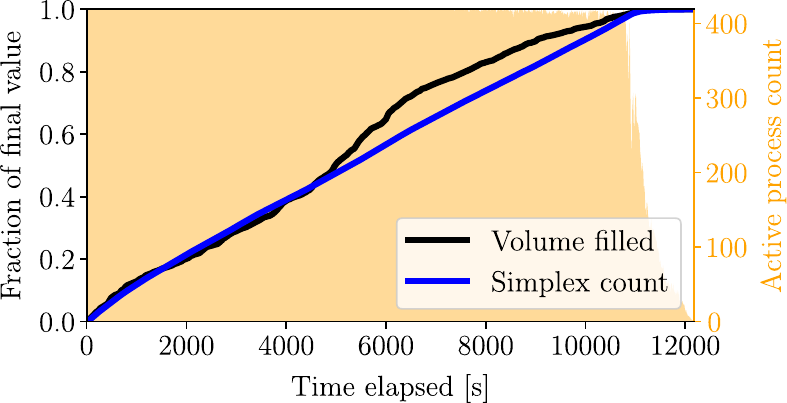}
  \caption{Algorithm~\ref{alg:lcss} progress plot. The partition complexity and
    the volume fraction of $\Theta$ comprised by completed tree branches
    steadily increase. The partitioning becomes serialized near the end.}
  \label{fig:progress}
\end{figure}

Table~\ref{tab:results} summarizes the output of Algorithm~\ref{alg:lcss} using
a sequence of increasingly tight $\epsilon$-suboptimality settings. We compute
$\epsilon_a$ as the largest cost among the vertices of a shrunk $\Theta$, i.e.
$\epsilon_{\mathrm{a}}=\max\{V^*(\theta)\mid\theta\in\mathcal V(s_{\mathrm{a}}\Theta)\}$. The
smaller $s_{\mathrm{a}}$ is, the more dense the partition will be in a
neighborhood of the origin. Figure~\ref{fig:progress} shows the progress of
Algorithm~\ref{alg:lcss} for the last row of Table~\ref{tab:results}. We can see
that the progress is mostly linear. Because multiple cores do not participate in
evaluating
lines~\ref{alg:lcss:line:check_suboptimality}-\ref{alg:lcss:line:parwhile_end}
for the same node $(\mathcal R,\delta)$, progress slows down near the end when
only a few nodes are left.

Figure~\ref{fig:depth_complexity} confirms that $\tau=\mathcal O(\log(\psi))$
using the proxy
$\psi=(\epsilon_{\mathrm{a}}/\bar\epsilon_{\mathrm{a}}+ \epsilon_{\mathrm{r}}/\bar\epsilon_{\mathrm{r}})^{-1}$,
where $\bar\epsilon_{\mathrm{a}}$ and $\bar\epsilon_{\mathrm{r}}$ are the
largest of the tested values. The leaf count is exponential in $\psi$. It
follows that $T_{\mathrm{cpu}}$ and $M$ are also exponential in $\psi$. Note
that the exponential increase in $T_{\mathrm{cpu}}$ can be offset by an
exponential increase in parallel core count, until a certain limit. Thus, our
method allows for reasonable $T_{\mathrm{wall}}$ runtimes.

\begin{figure}[t]
  \centering
  \begin{subfigure}[t]{1\columnwidth}
    \centering
    \includegraphics[width=1\linewidth]{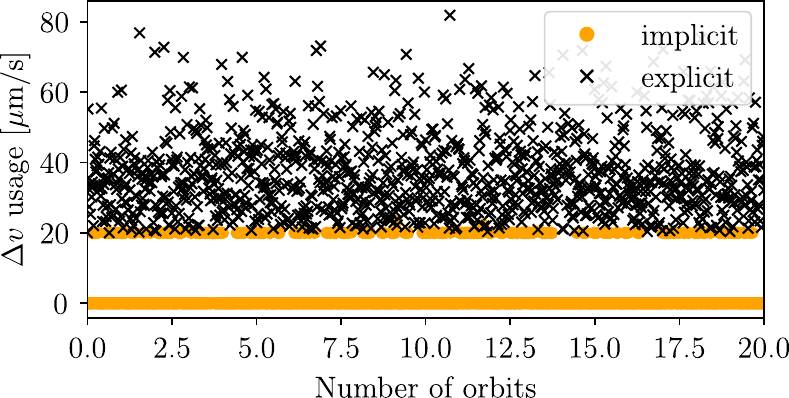}
    \caption{Input 2-norm history for \texttt{cwh\_z} with $s_{\mathrm{a}}=0.5$
      and $\epsilon_{\mathrm{r}}=2$.}
    \label{fig:input_2norm_coarse}
  \end{subfigure}%
  \hfill%
  \begin{subfigure}[t]{1\columnwidth}
    \centering
    \includegraphics[width=1\linewidth]{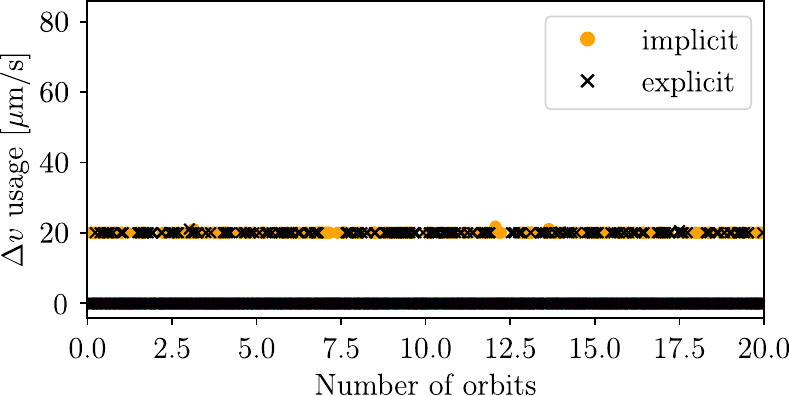}
    \caption{Input 2-norm history for \texttt{cwh\_z} with $s_{\mathrm{a}}=0.01$
      and $\epsilon_{\mathrm{r}}=0.01$.}
    \label{fig:input_2norm_refined}
  \end{subfigure}%
  \caption{Comparison of control input histories for a coarse and a refined
    $\epsilon$-suboptimal partition. By reducing $\epsilon_{\mathrm{a}}$ and
    $\epsilon_{\mathrm{r}}$, explicit MPC approaches the behavior of implicit
    MPC.}
  \label{fig:input_2norm}
\end{figure}

\begin{figure*}[t]
  \centering
  \begin{subfigure}[t]{.48\textwidth}
    \centering
    \includegraphics[width=1\linewidth]{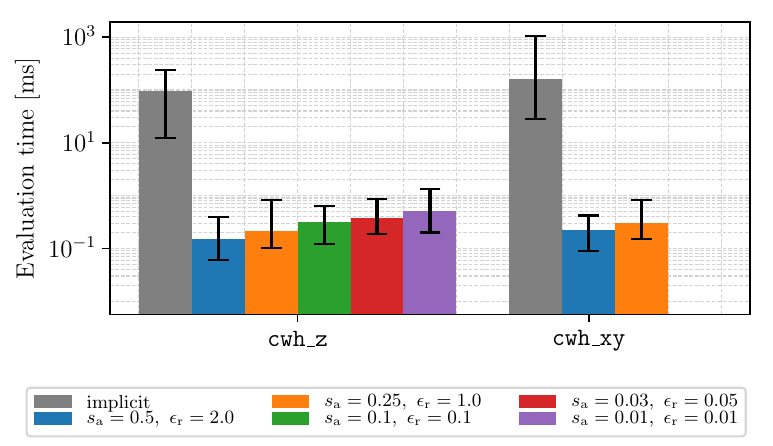}
    \caption{MPC on-line evaluation time. Bars show the mean while error bars
      shown the minimum and maximum values.}
    \label{fig:evaltime}
  \end{subfigure}%
  \hfill%
  \begin{subfigure}[t]{.48\textwidth}
    \centering
    \includegraphics[width=1\linewidth]{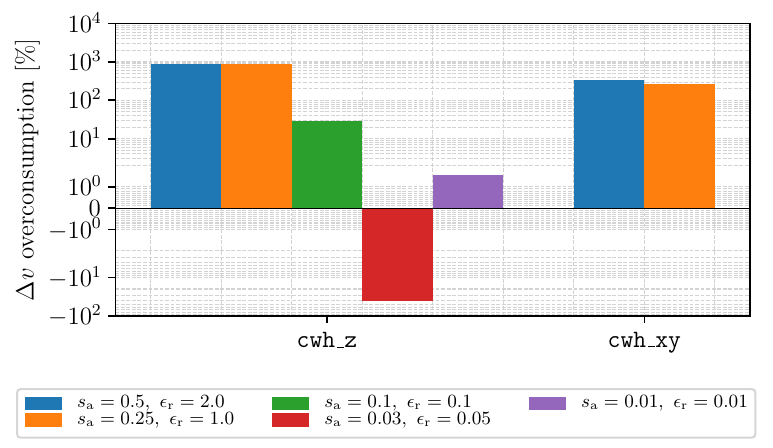}
    \caption{Overconsumption of fuel with respect to implicit MPC due to
      $\epsilon$-suboptimality. Implicit MPC uses $\approx 4$~mm/s over 20
      orbits.}
    \label{fig:fuel}
  \end{subfigure}%
  \caption{Comparison of the proposed semi-explicit and explicit implementations
    to implicit MPC in terms of (\protect\subref{fig:evaltime}) on-line control
    input computation time and (\protect\subref{fig:fuel}) total fuel
    consumption over $20$ orbits.}
\end{figure*}

Figure~\ref{fig:evaltime} show statistics for the on-line control input
computation time. For implict MPC, this is the mixed-integer solver time. For
explicit MPC, it is the time to evaluate~\eqref{eq:explicit_implementation},
which involves querying the partition tree. Statistics are computed by uniformly
randomly sampling 1000~values of $\theta\in\Theta$. As expected from
Theorem~\ref{theorem:onlineevalcomplexity}, explicit MPC can be up to three
orders of magnitude faster. Importantly, explicit MPC provides a real-time
guarantee given by the time that it takes to traverse the tree to the deepest
leaf. The implicit approach may be arbitrarily slower for some values of
$\theta$, subject to the success of the mixed-integer solver's heuristics.

Figure~\ref{fig:fuel} quantifies the fuel consumption suboptimality with respect
to implicit MPC. The data is collected based on a 20~orbit simulation where the
satellite is initialized at the origin, i.e. with zero control error. As
expected, partitions with a tighter $\epsilon$-suboptimality setting perform
better. Importantly, explicit MPC can \textit{outperform} implicit MPC since the
control scheme is finite horizon while fuel is an integrated quantity. This is
the case for \verb|cwh_z| with $s_{\mathrm{a}}=0.03$ and
$\epsilon_{\mathrm{r}}=0.05$, which achieves $\approx 40~\%$ fuel reduction. The
source of this reduced fuel consumption is clearly visible in
Figure~\ref{fig:input_2norm}, where one can see that with a tighter
$\epsilon$-suboptimality setting, explicit MPC better reproduces the optimal
behavior of implicit MPC.

\subsection{Improving the Storage Memory Requirement}
\label{subsec:improving_storage_memory_requirement}

Our implementation is not optimized for storage size, so the $M$ values in
Table~\ref{tab:results} are far greater than necessary. A more efficient storage
model is as follows. Given a simplex $\mathcal R\subset\reals^p$ and its
vertices $v_i\in\mathcal V(\mathcal R)$, a parameter $\theta\in\mathcal R$ if an
only if $\alpha = H_{\mathcal R}^{-1}(x-v_1)\ge 0$, $\bm{1}^\transp\alpha\le 1$,
where the $i$-th column of $H_{\mathcal R}$ equals $v_{i+1}-v_1$. Since $\Theta$
and hence $\mathcal R$ are full-dimensional, $H_{\mathcal R}$ is invertible. By
leveraging mutual exclusivity of the partition cells, we can thus store a matrix
$H\in\reals^{p\times p}$ and a vector $v_1\in\reals^p$ for each ``left'' child
node. For each leaf, in order to evaluate \eqref{eq:explicit_implementation} we
store the $p+1$ $\epsilon$-suboptimal decision vectors $x_j^*$ at its
vertices. Assuming a perfect binary tree and that $\Theta$ is a simplex, the
improved storage size is:
\begin{equation}
  \label{eq:optimized_storage_size}
  M^*\approx\frac{3}{2}\lambda\mu_{\mathrm{f}}p(p+1)+\lambda(p+1)\hat n\mu_{\mathrm{f}},
\end{equation}
where $\mu_{\mathrm{f}}$ is the floating point size and $\hat n\le n$ is the
dimension of the part of the decision vector that is necessary to compute the
control input (i.e. the first control input for MPC). For the examples in
Table~\ref{tab:results}, $M^*$ is 3 to 10 times less than $M$. Greater economy
is possible by eliminating further redundancy in the stored vertices and optimal
decision vectors.

\section{Conclusion and Future Work}
\label{sec:conclusion}

This paper presented a partitioning algorithm for generating explicit solutions
of a very general class of multiparametric mixed-integer convex programs to
within a given suboptimality tolerance. We showed that the positivity of a novel
cost function overlap metric is necessary and sufficient for algorithm
convergence. To the best of our knowledge, this is the first deep theoretical
insight into the fundamental driver of convergence rate and partition complexity
of suboptimal explicit MPC. In future work it will be interesting to prove the
stability of the resulting control law along the lines of
\cite{MunozDeLaPena2004} and \cite{MunozDeLaPena2006}, and to see if the
selection of $\epsilon_{\mathrm{a}}$ and $\epsilon_{\mathrm{r}}$ could be
automated to ensure convergence.

\section{Acknowledgments}

This research was partially supported by the National Science Foundation
(CMMI-1613235). The use of advanced computational, storage, and networking
infrastructure was provided by the Hyak supercomputer system and funded by the
STF at the University of Washington. The authors would like to thank Martin
Cacan, David S. Bayard, Daniel P. Scharf, Jack Aldrich and Carl Seubert of the
NASA Jet Propulsion Laboratory, California Institute of Technology, for their
helpful insights and discussions.

\bibliographystyle{ieeetr}
\bibliography{bibliography.bib}

\end{document}